\documentclass[10pt]{amsart}
\usepackage{amsmath,amssymb,amscd,latexsym}

\theoremstyle{plain}
\newtheorem{Thm}{Theorem}
\newtheorem{Cor}{Corollary}
\newtheorem{Prop}{Proposition}[section]
\newtheorem{Lem}[Prop]{Lemma}
\newtheorem{Def}{Definition}

\input{epsf}
\theoremstyle{definition}

\begin{document}

%\begin{document}
%Topmatter
\title
{The Kauffman skein module of the connected sum of $3$-manifolds}

\author{Jianyuan K. Zhong}
\address{Department of Mathematics \& Statistics\\
   Louisiana Tech University\\
    Ruston, LA 71272}
\email{kzhong@coes.latech.edu}
%\subjclass[2000]{57M}
%\date{June 11, 2001}
\date{September 17, 2002}
%End topmatter
\begin{abstract}
Let $k$ be an integral domain containing the invertible elements
$\alpha $, $s$ and $\frac{1}{s-s^{-1}}$.
If $M$ is an oriented $3$-manifold,  let $K(M)$ denote the
Kauffman skein module of $M$ over $k$. Based on the work on Birman-Murakami-Wenzl algebra by Beliakova
and Blanchet \cite{BB2000}, we give an ``idempotent-like'' basis for the Kauffman skein module of handlebodies. Gilmer and Zhong \cite{GZ01} have studied the Homflypt skein modules of a connected sum of two $3$-manifolds, here we study the case for the Kauffman skein module and show that $K(M_1\# M_2)$ is
isomorphic to $K(M_1)\otimes K(M_2)$ over a certain localized
ring, where $M_1\# M_2$ is the connected sum of two manifolds $M_1$ and $M_2$.
\end{abstract}

\maketitle
\section{Introduction}
%Let $k={\mathbb Z}[\alpha^{\pm 1}, s^{\pm 1}]$ containing $\frac{1}{s-s^{-1}}$.
Let $k$ be an integral domain containing the invertible elements
$\alpha $ and $s$. Moreover, we assume $s-s^{-1}$ is invertible in $k$.

We will work with unoriented framed links. An unoriented link is one equipped with a nonzero normal vector field up to
homotopy. The links described by figures in this paper are
assigned the vertical framing which points towards the reader.

Let $M$ be a smooth, compact and oriented $3$-manifold.
\begin{Def}
{\bf The Kauffman skein module} of $M$, denoted by $K(M)$, is the
$k$-module freely generated by isotopy classes of framed links in
$M$, including the empty link, modulo the following Kauffman skein relations:
$$(i)\quad \raisebox{-3mm}{\epsfxsize.3in\epsffile{left1.ai}}\quad
        -\quad\raisebox{-3mm}{\epsfxsize.3in\epsffile{right1.ai}}\quad
=\quad (\ s - \
s^{-1})\Bigl(\quad\raisebox{-3mm}{\epsfxsize.3in\epsffile{parra1.ai}}-\quad\raisebox{-3mm}{\epsfxsize.3in\epsffile{parra2.ai}}\quad\Bigr)\quad
,$$
$$(ii)\quad \raisebox{-3mm}{\epsfxsize.35in\epsffile{framel1.ai}}
        \quad=\quad \alpha \quad \raisebox{-3mm}{\epsfysize.3in\epsffile{orline1.ai}}\quad  ,
$$
$$(iii)\quad L\ \sqcup
\raisebox{-2mm}{\epsfysize.3in\epsffile{unknot1.ai}}\quad =\quad
\delta\quad L\quad,
$$
where $\delta=({\dfrac{\alpha-\alpha^{-1}} {\ s - \ s^{-1}}}+1)$. The last relation follows from the first two when $L$ is nonempty.
\end{Def}
\noindent Equations $(i)$ and $(ii)$ are local relations, diagrams in equation $(i)$ are identical except where shown, so are diagrams in equation $(ii)$. A trivial closed curve in $(iii)$ is a curve
which is nul-homotopic and contains no crossing. 

We list here some important basic facts on Kauffman skein modules. (1) An embedding $f: M\to N$ of $3$-manifolds
induces a well defined homomorphism $f_{*}: K(M)\to  K(N)$. (2) If we obtain
$N$ by adding a $3$-handle to $M$, then the embedding $i:
M\to N$ induces an isomorphism $i_{*}:  K(M)\to  K(N)$. (3) If we add a $2$-handle to $M$ to obtain $N$, then the embedding $i: M\to
N$ induces an epimorphism $i_{*}:  K(M)\to  K(N)$. (4) $ K(M_1 \sqcup M_2)\cong  K(M_1)\otimes  K(M_2)$, where $M_1
\sqcup M_2$ is the disjoint union of two $3$-manifolds $M_1$ and
$M_2$.

Let $\lambda$ be a Young diagram associated to a partition of $n$, where $\lambda=(\lambda_1\geq \dots \lambda_p \geq 1)$ with $\lambda_{1}+\dots +\lambda_{p }=n$, then $\lambda$ has $|\lambda|=n$ cells indexed by \{$(i, j), \ 1\leq i\leq p$,
$1\leq j\leq \lambda_i$\}; if $c$ is the cell of index $(i, j)$ in $\lambda$, its content
$cn(c)$ is defined by $cn(c)=j-i.$ We 
define a scalar $c_\lambda$ associated to $\lambda$ as,
$$c_\lambda=(\ s - \
s^{-1})(\alpha\sum_{c\in \lambda}s^{2cn(c)}-\alpha^{-1}\sum_{c\in
\lambda}s^{-2cn(c)}).$$

Let ${\mathcal{I}}$ denote the submonoid  of the multiplicative monoid of 
${\mathbb{Z}}[\alpha,s]$ generated by $\alpha$, $s$, ${s^{2n}-1}$ and $1\pm\alpha s^n$ for all integers $n>0,$ and $c_{\lambda}$ for any Young diagrams $\lambda$.  Let ${\mathcal R}$ be  ${\mathbb{Z}}[\alpha,s]$ localized at  ${\mathcal I}$ \cite[7.2]{J}

\begin{Thm} Over $k={\mathcal {R}}$, 
 $$K(M_1 \# M_2)\cong K(M_1)\otimes K(M_2).$$
\end{Thm}
\noindent {\em Remark}: Przytycki has proved the analog of this
result for the Kauffman bracket skein module \cite{JP}, and Gilmer
and Zhong have settled the case for the Homflypt skein module
\cite{GZ01}. We follow a similar general outline as in \cite{GZ01} but with a different approach using the Kauffman skein module of the handlebodies.

\begin{Cor} Let $k={\mathbb Q}(\alpha, s)$, the field
of rational functions in $\alpha, s$, then
 $$K(M_1 \# M_2)\cong K(M_1)\otimes K(M_2).$$
\end{Cor}

From \cite{ZB01}, we know that over $k={\mathbb Q}(\alpha, s)$, the Kauffman skein module of $S^1\times S^2$ is freely generated by the empty link. Let $\#^m S^1\times S^2$ denote the connected sum of $m$ copies of $S^1\times S^2$, we have the following corollary.
\begin{Cor} When $k={\mathcal {R}}$ or $k={\mathbb Q}(\alpha, s)$, $K( \#^m S^1\times S^2)$ is a free $k$-module generated by the empty link.
\end{Cor}

\begin{Def}
{\bf The relative Kauffman skein module}. Let $X=\{x_1,
x_2,\cdots, x_{n},$ $y_1, y_2,\cdots, y_{n}\}$ be a finite set of $2n$
framed points in the boundary $\partial M$ of $M$. Define the relative
skein module $K(M, X)$ to be the $k$-module freely generated by relative
framed links in $(M,\partial M)$ such that $L \cap
\partial M=\partial L=\{x_i, y_i\}$ with the induced framing,
considered up to an ambient isotopy fixing $\partial M$ modulo the
Kauffman skein relations.
\end{Def}

We have the following version of Theorem 1 for relative Kauffman skein modules. In this case we too have to work over ${\mathbb{Q}}(\alpha,s)$, the field of fractions of ${\mathbb{Z}}[\alpha^{\pm 1},s^{\pm 1}]$, because we do not know whether the relative skein module of a handlebody is free. We conjecture that it is also free.

\begin{Thm}
Let $M_1$ and $M_2$ be connected oriented 3-manifolds, and let  $X=\{x_1, x_2,\cdots, x_{n}\}$ and $Y=\{y_1, y_2,\cdots, y_{n}\}$ be two sets of framed points in $\partial M_1$ and $\partial M_2$, respectively, then
$$K(M_1 \# M_2, X\cup Y)\cong K(M_1, X)\otimes K(M_2,Y).$$
\end{Thm}

In section 2, we give an ``idempotent-like'' basis for a handleboy of genus $g$ using the relative Kauffman skein module of a disc and the work of Beliakova and Blanchet \cite{BB2000} on the Birman-Murakami-Wenzl algebras. In section 3, we prove Theorem 1 first in the case of two handlebodies, then in the general case. In section 4, we give an example in the connected sum of two copies of $S^1\times S^2$. In section 5, we prove a lemma which supports the Epimorphism lemma 3.2 in section 3. The last three sections follow similar outlines of \cite{GZ01}.

\section{The Kauffman skein module of handlebodies}

\subsection {The Birman-Murakami-Wenzl category and the
Birman-Murakami-Wenzl algebra $K_n$}

In {\it the Birman-Murakami-Wenzl category $K$} \cite{BB2000}, an
object in $K$ is a disc $D^2$ equipped with a finite set of framed points; if $\beta=(D^2, l_0)$ and
$\gamma=(D^2, l_1)$ are objects, the module $Hom_K(\beta, \gamma)$
is the relative Kauffman skein module $K(D^2\times [0,1], l_0\times 0 \amalg l_1\times 1)$. If $\lambda$ is a
Young diagram, we denote by $\square_\lambda$ the object
of the category $K$ formed with one point assigned for each cell of
$\lambda$. We will use the notation $K(\beta, \gamma)$ for $Hom_K(\beta,
\gamma)$. For composition, it's by stacking $f$ on
the top of $g$.
$$K(\beta, \gamma)\times K(\gamma, \delta)\to K(\beta, \delta),$$
$$(f,g)\to fg$$
Note that Beliakova and Blanchet choose to stack the second one on
the top of the first. Here we follow the convention as in
\cite{AM98}.

As a special case, in the cylinder $D^2\times I$, If $X_n$ is a
set of $n$ distinct framed points on a diameter $D^2\times \{1\}$
and $Y_n$ is a set of $n$ distinct framed points on a diameter
$D^2\times \{0\}$, then the relative Kauffman skein module
$K(D^2\times I, X_n\amalg Y_n)$ is isomorphic to the {\it
Birman-Murakami-Wenzl algebra} $K_n$, $K_n$ is generated by the
identity ${\mathbf{1}}_n$, positive transpositions $e_1, e_2,
\cdots, e_{n-1}$ and hooks $h_1, h_2, \cdots, h_{n-1}$ as the
following:
$$e_i=\quad \raisebox{-3.5mm}{\epsfxsize.0in\epsffile{ei.ai}}$$
$$h_i=\quad \raisebox{-3.5mm}{\epsfxsize.0in\epsffile{hi.ai}}$$
for $1\leq i\leq n-1$.
 
Let $I_n=\{(a\otimes {\mathbf 1}_1)h_{n-1}(b\otimes {\mathbf 1}_1):\ a,
b\in K_{n-1}\}$ be the ideal generated by $h_{n-1}$.

We recall that in the case of the Homflypt skein module, the
$n$th Hecke algebra of type $A$ is realized geometrically \cite{AM98} by the
relative Homflypt skein module of the cylinder $D^2\times I$ with a set of $n$ distinct input framed points on a
diameter $D^2\times \{1\}$ and a set of $n$ distinct
output framed points on a diameter $D^2\times \{0\}$. 

It has been revealed \cite{BB2000} that the algebras $K_n$ and $H_n$ are closely related. On the one hand, modulo the ideal $I_n$, there is a canonical projection map from $K_n$ to $H_n$:
$$\pi_n: K_n\to H_n.$$ 
On the other hand, we have the following theorem \cite{BB2000}. 

\begin{Thm}
There is a multiplicative homomorphism $s_n:H_n\to K_n$, such
that
$$\pi_n\circ s_n=id_{H_n},$$
$$s_n(x)y=ys_n(x)=0, \ \forall x\in H_n, \forall y \in I_n.$$
\end{Thm}

i.e. $K_n=s_n(H_n)\oplus I_n$.

This theorem was proved over ${\mathbb Q}(\alpha, s)$, the field
of rational functions in $\alpha, s$. See details of the proof in
\cite{BB2000}. But we comment that this theorem holds if elements of the form $1-s^{2n}$ and $1\pm\alpha s^n$, for all integers $n>0$,  are invertible in $k={\mathbb Z}[\alpha^{\pm 1}, s^{\pm 1}]$.

Let $\lambda$ be a Young diagram of size $n$. In the
Hecke category $H_{\square_{\lambda}}$, we have a Young idempotent
$y_{\lambda}$, whose definition and properties are in
\cite[Chapter 3]{GZ00} \cite{cB98}. Let $y_{\lambda}^{*}$ be the
flattened version of $y_{\lambda}$ according to Aiston and Morton \cite{AM98}, who gave an example of a three-dimensional view of the quisi-idempotent associated with the Young diagram $\lambda$, from which the idempotent $y_{\lambda}$ is obtained after normalization. They assigned $\lambda$ a standard tableau where the cells of $\lambda$ are numbered from $1$ to $|\lambda|$ along the rows, and a similar standard tableau is assigned for the transpose of $\lambda$ (whose rows are the columns of $\lambda$). These two tableaux induce a permutation $\pi_\lambda$ between the cells in the rows of $\lambda$ and in the columns of the transpose of $\lambda$. A flattened version of the idempotent is obtained by sliding the rows apart on the top of the three-dimensional diagram and sliding the columns apart at the bottom, according to the ordering of the strings induced by $\pi_\lambda$. (See details in \cite{AM98}.)

We denote the homomorphic image of $y_{\lambda}^{*}$ under $s_n$ by $\tilde{y}_\lambda$, $\tilde{y}_\lambda=s_n(y_{\lambda}^{*})$. We denote the natural closure of  $\tilde{y}_\lambda$ in $S^3$
by $\widehat{\widetilde{y}_\lambda}$.

\subsection{A basis for the Birman-Murakami-Wenzl algebra
$K_n$}

Here we mainly summarize some results by Beliakava and Blanchet
\cite{BB2000}.

An up and down tableau of length $n$ and
shape $\Lambda_n$ is a sequence $\Lambda=(\Lambda_1, \cdots, \Lambda_n)$ of Young
diagrams with every two consecutive Young diagrams $\Lambda_i$
and $\Lambda_{i+1}$ distinct by exactly one cell. Note that the size of $\Lambda_n$ is either $n$ or less than $n$
by an even number.

For an up and down tableau $\Lambda$ of length $n$, let
$\Lambda'$ be the tableau of length $n-1$ obtained by removing the
last Young diagram in the sequence $\Lambda$. Beliakava and Blanchet inductively define
$a_{\Lambda}\in K(n, \square_{\Lambda})$ and $b_{\Lambda}\in K(
\square_{\Lambda}, n)$ by
$$a_1=b_1={\mathbf{1}}_1.$$
If $|\Lambda_n|=|\Lambda_{n-1}|+1$, then
$$a_{\Lambda}=(a_{\Lambda'}\otimes
{\mathbf{1}}_1){\widetilde{y}}_{\Lambda_n},$$
$$b_{\Lambda}={\widetilde{y}}_{\Lambda_n}(b_{\Lambda'}\otimes
{\mathbf{1}}_1);$$ if $|\Lambda_n|=|\Lambda_{n-1}|-1$, then
$$a_{\Lambda}=\frac{<\Lambda_n>}{<\Lambda_{n-1}>}(a_{\Lambda'}\otimes
{\mathbf{1}}_1)({\widetilde{y}}_{\Lambda_n}\otimes \cup),$$
\hspace{3.5cm} $b_{\Lambda}=({\widetilde{y}}_{\Lambda_n}\otimes
\cap)(b_{\Lambda'}\otimes {\mathbf{1}}_1).$

Here $<\lambda>$ is the Kauffman polynomial of
$\widehat{\widetilde{y_\lambda}}$ in $S^3$ called the quantum dimension associated with
$\lambda$. Note $<\lambda>$ is invertible if elements of the form $1-s^{2n}$ and $1\pm\alpha s^n$, for all integers $n>0$,  are invertible in $k={\mathbb Z}[\alpha^{\pm 1}, s^{\pm 1}]$ \cite{BB2000}.

\begin{Thm}
A basis for $K_n$ is given by the family $a_{\Lambda}b_{\Gamma}$ for all up and down tableaux
$\Lambda, \Gamma$ of length $n$, such that $\Lambda_n=\Gamma_n$.
\end{Thm}

See details of the proof in \cite{BB2000}. From the
corresponding properties in the Hecke category, it follows that if $\Lambda=\Gamma$ ( $\Lambda_i=\Gamma_i$ for $1\leq i \leq
n$), then $b_{\Gamma}a_{\Lambda}=\widetilde{y}_{\Lambda_n}$;
otherwise $b_{\Lambda}a_{\Gamma}=0$.

Let $\Lambda=(\Lambda_1, \cdots, \Lambda_n), \
\Gamma=(\Gamma_1,\cdots, \Gamma_n)$ be two up and down tableaux of
length $n$ with $\Lambda_n=\Gamma_n$, we give the following figures for illustration.

(1) if $|\Lambda_n|=|\Lambda_{n-1}|+1$ (so is
$|\Gamma_n|=|\Gamma_{n-1}|+1$), then
$$a_{\Lambda}b_{\Gamma}=\raisebox{-30mm}{\epsfxsize.0in
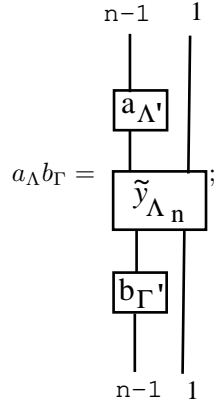
\epsffile{alabth.ai}};$$
(2) if $|\Lambda_n|=|\Lambda_{n-1}|-1$ (so is
$|\Gamma_n|=|\Gamma_{n-1}|-1$), then
$$a_{\Lambda}b_{\Gamma}=\frac{<\Lambda_n>}{<\Lambda_{n-1}>}\quad
\raisebox{-35mm}{\epsfxsize.0in
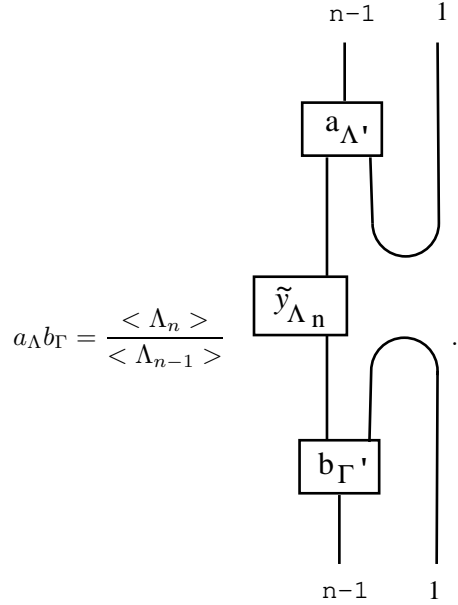
\epsffile{alabth1.ai}}.$$

\subsection{The Kauffman skein module of handlebodies}
\begin{Lem}
Let $A$, $B$, $C$ be three sets of framed points on the boundary
of $D^2$ as shown, where $A$, $B$, $C$ contain the numbers of $a$, $b$, $c$
framed points, respectively. Assume $a+b+c$ is a nonnegative even
integer, then the relative Kauffman skein module $K(D^2, A\cup
B\cup C)$ has a generating set with generators of the form:
$$\raisebox{-3mm}{\epsfxsize.0in
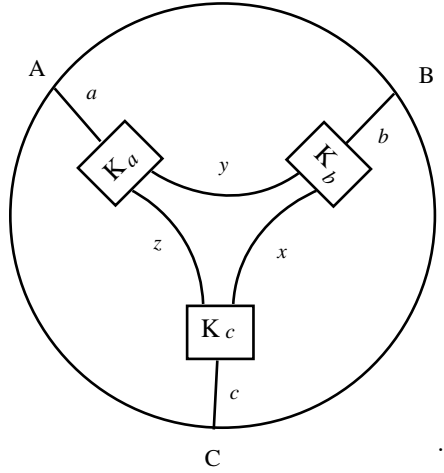
\epsffile{d3.ai}}.$$
\end{Lem}
\noindent where $\raisebox{-3mm}{\epsfxsize.0in\epsffile{kd.ai}}$ represents 
all possible basis elements of the $i$th Birman-Murakami-Wenzl algebra
$K_i$ given in Theorem 4, and $a+b+c=x+y+z$. In the diagram, a string associated with a nonnegative integer $d$ denotes $d$ parallel strings. When $a+b+c=0$,
$K(D^2, A\cup B\cup C)$ is generated by the empty diagram.

{\bf{Sketch proof.}}
Given a framed unoriented link $L$ in $K(D^2, A\cup B\cup C)$, we can first apply the Kauffman skein relations to rewrite $L$ as a linear combinations of elements in $K(D^2, A\cup B\cup C)$ without closed components; then by an induction on the number of crossings in the diagram and rearranging strings coming from the same set of boundary points into the same box, eventually, we can rewrite $L$ as a linear sum of elements in the form given in the lemma.

Now we define a triple of nonnegative integers $a, b, c$ to be an admissible triple if $a+b+c$ is a nonnegative even integer. We will use a box with an integer $i$ next to the box to represent $\raisebox{-3mm}{\epsfxsize.0in\epsffile{kd.ai}}$. We introduce the following simpler diagram to represent the figure given in the previous lemma:
$$\raisebox{-3mm}{\epsfxsize.0in
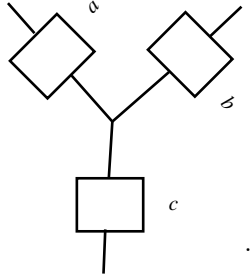
\epsffile{d31.ai}}.$$

We call the above figure an admissible triad. Note this figure represents a group of elements with each box being a possible basis element from the corresponding $K_i$, for $i=a,b,c$.

Let ${\mathcal H }_g$ be a handlebody of genus $g$, ${\mathcal H }_g$ can be defined as $F_g \times I$, where $F_g$ is a
closed disc deleting $g$ open discs and $I$ is the unit
interval$[0,1]$. Consider a trinion (or called pair of pants) decomposition of
${\mathcal{H}}_g$, we can paste the
generators of $K(D^2, A\cup B\cup C)$ from the previous lemma through the trinion decomposition of ${\mathcal H }_g$ and obtain a
generating set for $K({\mathcal{H}}_g)$.

Let ${\mathcal G}$ be the set of elements in $K({\mathcal{H}}_g)$ of the following form,
$$\raisebox{-4mm}{\epsfxsize.0in
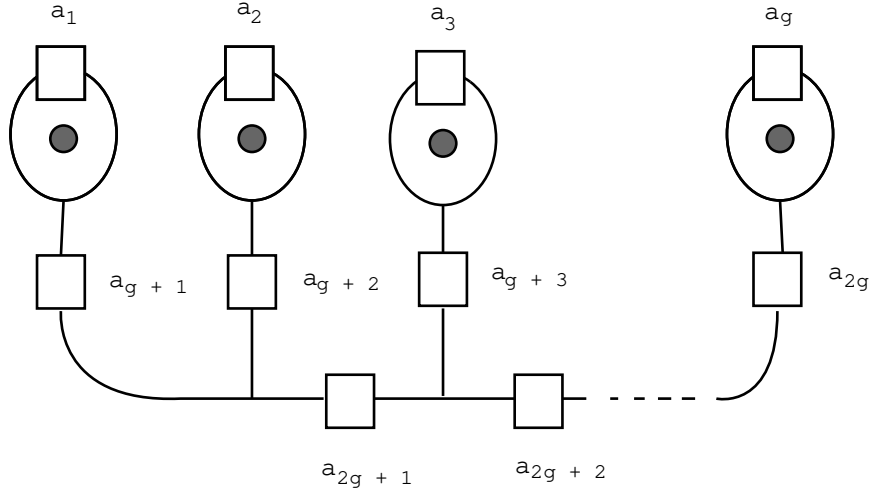
\epsffile{hgbasis.ai}}$$
where $a_1, a_2, \cdots, a_g$, $a_{g+1}, a_{g+2},\cdots, a_{2g}$, $a_{2g+1}, a_{2g+2}, \cdots, a_{3g-2}$ are nonnegative integers, which form all admissible triads in the picture. A gray circle in the figure indicates a deleted open disk.

\begin{Thm}
$K({\mathcal{H}}_g)$ has a generating set given by  ${\mathcal G}$.
\end{Thm}

\begin{Cor} Over ${\mathbb Q}(\alpha, s)$,
$K({\mathcal{H}}_g)$ has a free basis given by a subset of  ${\mathcal G}$.
\end{Cor}
\begin{proof}
Note that $K({\mathcal{H}}_g)$ is a free module over ${\mathbb Q}(\alpha, s)$. By Zorn's Lemma, there is a maximum linearly independent subset of ${\mathcal G}$ containing the empty framed link, which generates $K({\mathcal{H}}_g)$ by a reductio ad absurdum proof. It follows that the maximal linearly independent generating set obtained above is a free basis.
\end{proof}

\section{Isomorphism for Handlebodies}

\subsection{Epimorphism for Handlebodies}

Let $D$ be a separating meridian disc of ${\mathcal{H}}_g$, let
$\gamma=\partial{D}$. 
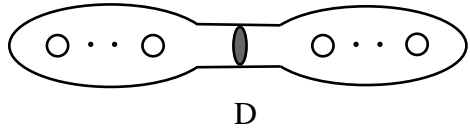
$$\raisebox{-3mm}{\epsfxsize.0in\epsffile{hm.ai}}$$

Let $V_D=[-1, 1]\times D$ denote a regular neighborhood of $D$ in
${\mathcal{H}}_g$, $V_D$ may be projected into a disc $D_p=[-1,
1]\times [0, 1]$.

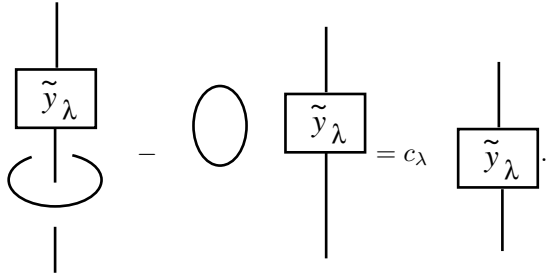
\begin{Lem}
$$\raisebox{-15mm}{\epsfxsize.0in\epsffile{alasab.ai}}\quad-\quad\raisebox{-13mm}{\epsfxsize.0in\epsffile{alasab0.ai}}=c_\lambda \quad
\raisebox{-12mm}{\epsfxsize.0in\epsffile{alasa.ai}}.$$
\end{Lem}

We will prove this lemma in section 5.

\begin{Lem}{\bf (The Epimorphism Lemma)} The embedding $i:\ {\mathcal{H}}_g-D \to ({\mathcal{H}}_g)_{\gamma}$ induces an epimorphism:
$$i_{*}: K({\mathcal{H}}_g-D)\twoheadrightarrow K(({\mathcal{H}}_g)_{\gamma}).$$
\end{Lem}
\begin{proof}
We will consider the following simple sliding relation from adding the $2$-handle along $\gamma$. In $V_D$, we illustrate the sliding relation as the following:
$$\raisebox{-13mm}{\epsfxsize.0in
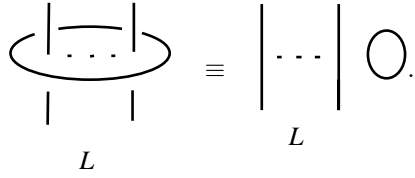
\epsffile{zno1.ai}}\quad\equiv\quad \raisebox{-10mm}{\epsfxsize.0in\epsffile{zno.ai}}.$$
(1) We first prove that any generator of $K({\mathcal{H}}_g)$ given by Theorem 5 which cuts $D$ nontrivially is zero in $K(({\mathcal{H}}_g)_{\gamma})$. Let $\beta$ be such a generator of $K({\mathcal{H}}_g)$, we assume that $\beta$ intersects $D$ nontrivially; then in $V_D$, $\beta$ has a basis element from $K_n$ for some $n>0$, which  contains an idempotent $\tilde{y}_\lambda$ for a nonempty Young diagram $\lambda$. Apply the above simple sliding relation to $\tilde{y}_\lambda$, from Lemma 3.1, we have $c_\lambda\beta=0$ in $K(({\mathcal{H}}_g)_{\gamma})$. As $c_\lambda$ is invertible, we conclude that $\beta=0$ in $K(({\mathcal{H}}_g)_{\gamma})$.
 
(2) In general, if $L$ is a framed link in $K({\mathcal{H}}_g)$, then $L$ can be written as a linear sum of the generators presented in Theorem 5. From (1), all generators cutting $D$ nontrivially in the linear sum of $L$ are zero, so $L\in K({\mathcal{H}}_g-D)$. The result follows.
\end{proof}

\subsection{Isomorphism for Handlebodies}

Let ${\mathcal{H}}_{g_1}$ and ${\mathcal{H}}_{g_2}$ be handlebodies with genera $g_1>0$, $g_2>0$ and let $g=g_1+g_2$. ${\mathcal{H}}_{g_1}\# {\mathcal{H}}_{g_2}$ is equal to ${\mathcal{H}}_{g_1+g_2}(={\mathcal{H}}_{g})$ with a 2-handle added along the boudary of the meridian disc $D$ separating ${\mathcal{H}}_{g_1}$ and ${\mathcal{H}}_{g_2}$. Let $\gamma=\partial D$, then ${\mathcal{H}}_{g_1}\# {\mathcal{H}}_{g_2}=({\mathcal{H}}_{g})_\gamma$. From the basic facts about the Kauffman skein modules, we have $K({\mathcal{H}}_{g_1})\otimes K({\mathcal{H}}_{g_2})\cong K({\mathcal{H}}_{g}-D)$. We will prove Theorem 1 by proving
$K({\mathcal{H}}_{g}-D)\cong K(({\mathcal{H}}_{g})_\gamma)$.

Let $I_0$ be the submodule of $K({\mathcal{H}}_g)$ given by the collection $\{p_{D} ( L ) - L\sqcup O \mid L \in K({\mathcal{H}}_g)\}$, where $O$ represents the unknot. Locally, we have the following diagram description.
$$p_{D}( L )=\raisebox{-12.5mm}{\epsfxsize.0in
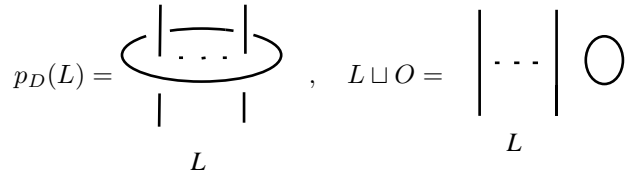
\epsffile{zno1.ai}}\quad, \quad L\sqcup O=\quad\raisebox{-10mm}{\epsfxsize.0in\epsffile{zno.ai}}$$

When $({\mathcal{H}}_g)_{\gamma}$ is obtained by adding a $2$-handle to ${\mathcal{H}}_{g}$ along $\gamma$, from a generalization of \cite{GZ00} section 2, we know that
$K(({\mathcal{H}}_g)_{\gamma})\cong K({\mathcal{H}}_{g})/ R$,
where $R$ is the submodule of $K({\mathcal{H}}_g)$ given by the
collection $\{\Phi '( z ) - \Phi ''( z ) \mid z \in
K({\mathcal{H}}_g, AB)\}$, where $A, \ B$ are two points on
$\gamma$, which decompose $\gamma$ into two intervals $\gamma '$
and $\gamma ''$, $z$ is any element of the relative skein module
$K({\mathcal{H}}_{g}, AB)$, and $\Phi '( z )$ and $\Phi ''( z )$
are given by capping off with $\gamma '$ and $\gamma ''$,
respectively, and pushed the resulting links back into
${\mathcal{H}}_{g}$.

\begin{Lem}
 $R=I_0$.
\end{Lem}

The proof of this lemma is very similar to that of Lemma 3.1 of \cite{GZ01}.

\begin{Cor}
The embedding ${\mathcal{H}}_{g} \to ({\mathcal{H}}_g)_{\gamma}$ induces an isomorphism
$$K({\mathcal{H}}_{g})/ I_0\cong K(({\mathcal{H}}_g)_{\gamma}).$$
\end{Cor}

Now to show that the embedding ${\mathcal{H}}_g-D \to ({\mathcal{H}}_g)_{\gamma}$ induces an isomorphism
$$K({\mathcal{H}}_g-D)\cong K(({\mathcal{H}}_g)_{\gamma}),$$
we only need to show $K({\mathcal{H}}_g-D) \cap I_0=0.$

\begin{Lem}
$$K({\mathcal{H}}_g-D) \cap I_0=0.$$
\end{Lem}

\begin{proof}
Lieberum \cite{Lie00} has shown that $K({\mathcal{H}}_g)$ is free, so $K({\mathcal{H}}_g-D)$ is also free. Let ${\mathcal F}={\mathbb Q}(\alpha, s)$, and let $K_{\mathcal F}(M)$ denote the Kauffman skein module of $M$ over ${\mathcal F}$. The map  $K({\mathcal{H}}_{g}-D) \rightarrow K_{\mathcal F}({\mathcal{H}}_{g}-D),$ induced by 
${\mathcal R}\rightarrow {\mathcal F}$ is injective. Let ${\mathcal {I}}_0=\{p_{D} ( L ) - L\sqcup O \mid L \in K_ {\mathcal {F}}({\mathcal{H}}_g)\}$. It is enough to show
 $K_{\mathcal F}({\mathcal{H}}_g-D) \cap {\mathcal I}_0=0.$

Let $\psi$ be the map from $K_{\mathcal F}({\mathcal{H}}_g) \rightarrow K_{\mathcal F}({\mathcal{H}}_g)$ given by $\psi(L)=p_D(L)-L\sqcup O$ for all $L\in K_{\mathcal F}({\mathcal{H}}_g)$.  
Note that $Image(\psi)={\mathcal I}_0$, and $\psi|_{K_{\mathcal F}({\mathcal{H}}_g-D)}=0$.

We make use of the free basis of $K_{\mathcal F}({\mathcal{H}}_g)$ obtained in Corollary 3. A similar argument to the proof of Lemma 3.1, we can show that each basis element of $K_{\mathcal F}({\mathcal{H}}_g)$ which cuts $D$ nontrivially is an eigenvector of $\psi$ with an invertible eigenvalue, hence $\psi$ is an injective map from $K_{\mathcal F}({\mathcal{H}}_g)/K_{\mathcal F}({\mathcal{H}}_g-D)$ to itself. Since $\psi|_{K_{\mathcal F}({\mathcal{H}}_g-D)}=0$, the result follows.
\end{proof}

\begin{Cor}
The embedding ${\mathcal{H}}_g-D \to ({\mathcal{H}}_g)_{\gamma}$ induces an isomorphism,
$$K({\mathcal{H}}_g-D)\cong K(({\mathcal{H}}_g)_{\gamma}).$$
\end{Cor}
This follows from the commutative diagram below from the above arguments,

$$\begin{CD}
K({\mathcal{H}}_{g}-D) @>>> K({\mathcal{H}}_{g}) 
\\
@V\text{epimorphism}VV      @V{I_0=R}VV
%@VVV      @VVV
\\
 K(({\mathcal{H}}_g)_{\gamma}) @>\cong>> K({\mathcal{H}}_{g})/I_0
\end{CD}
$$

\subsection{Isomorphism in the general case}

If $M$ is a closed connected oriented $3$-manifold, $M$ can be be obtained from a handlebody ${\mathcal H}$ by adding some $2$-handles and one 3-handle. As adding or removing 3-balls from the interior of a 3-manifold does not change its Kauffman skein module, we prove Theorem 1 assuming that $M_1$ and $M_2$ are connected 3-manifolds with boundary; in which case, each $M_i$ is obtained from the handlebody ${\mathcal H}_{g_i}$ by adding some $2$-handles. 
Let $g=g_1+g_2$. Let $N$ be the manifold obtained by adding both sets of  $2$-handles to the boundary connected sum of  ${\mathcal H}_{g_1}$ and ${\mathcal H}_{g_2}$ which may be identified as ${\mathcal H}_{g}$. Let $D$ be the disc in 
${\mathcal H}_{g}$ separating ${\mathcal H}_{g_1}$ from ${\mathcal H}_{g_2}.$  We let $\gamma=\partial{D},$ so ${\mathcal H}_{g_1}\# {\mathcal H}_{m_2}=({\mathcal H}_{g})_{\gamma}.$ Let  $N_{\gamma}$ denote the manifold obtained by adding a 2-handle to $N$ along a curve  $\gamma$ in $\partial N.$
We consider $M_1\# M_2$ as obtained from $({{\mathcal H}}_{g})_{\gamma}$ by adding those $2$-handles. Thus $N - D = M_1 \sqcup M_2,$  and $M_1\# M_2=N_{\gamma}.$
 
\begin{Thm}
The embedding  $N-D \rightarrow N_{\gamma}$ induces an isomorphism
$$K(N-D)\cong K(N_{\gamma}).$$
\end{Thm}

{\bf Remark:} The proof is an analogy to that of Theorem 5 in \cite{GZ01} by an induction on the number of the $2$-handles to be added to $({{\mathcal H}}_{g})_{\gamma}$ to obtain $N_{\gamma}.$ 

\begin{Cor} Let $B_1$ and $B_2$ be the 3-balls removed from $M_1$ and $M_2$ while forming the connected sum $M_1\# M_2.$ The embedding  $(M_1 - B_1)\sqcup (M_2 - B_2) \rightarrow M_1\# M_2$ induces an isomorphism
$$K(M_1)\otimes K(M_2)\cong K(M_1\# M_2).$$
\end{Cor}
\begin{proof}
Since $K(M-D) \cong K(M_1)\otimes K(M_2)$.
\end{proof}

Similarly, the above corollary holds whether or not $M_1$ or $M_2$ have boundary.

\section {An example in $S^1\times S^2\# S^1\times S^2$}

In \cite{ZB01}, we showed that $K(S^1\times S^2)$ is a free ${\mathcal R}$-module generated by the empty link. It follows that $K(S^1\times S^2\# S^1\times S^2)$ is also a free module generated by the empty link.
Let $K_{knot}$ be a knot in $S^1\times S^2\# S^1\times S^2$ given by the following diagram,
$$\raisebox{-3mm}{\epsfxsize.0in
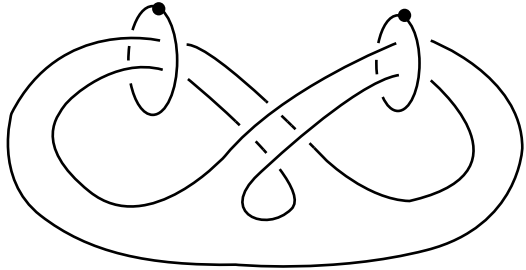
\epsffile{exs1s2.ai}}$$
Here the two circles each with a dot are a framed link description of $S^1\times S^2\# S^1\times S^2$. 
Note this same knot was studied with respect to the Kauffman Bracket skein modules in 
\cite{pG98} and Homflypt skein module in \cite{GZ01}.

In $K(S^2\times S^1 \# D^3, 4pts)$, isotopy  yields,
\begin{equation}
\raisebox{-12mm}{\epsfxsize.0in\epsffile{ex1.ai}}= \raisebox{-12mm}{\epsfxsize.0in\epsffile{ex2.ai}}\tag{I}
\end{equation}

Let $\lambda$ be the Young diagram of size two with the two cells on one row,
let $\mu$ be the Young diagram of size two with the two cells on one column.
$$\lambda=\quad\raisebox{-3mm}{\epsfysize0in\epsffile{2l.ai}}\quad, \quad\quad\mu=\quad\raisebox{-5mm}{\epsfysize0in\epsffile{2b.ai}}.$$
Recall in the Hecke algebra $H_2$, the idempotents corresponding
to  $\lambda$ and $\mu$ are the symmetrizer $f_2$ and antisymmetrizer $g_2$
given by Beliakova and Blanchet
\cite{BB2000}[Chapter 2] with
$$f_2=\frac{s^{-1}}{s+s^{-1}}\raisebox{-3mm}{\epsfxsize.3in\epsffile{parra.ai}}\quad + \frac{s^{-1}}{s+s^{-1}}\quad \raisebox{-3mm}{\epsfxsize.3in\epsffile{left.ai}}$$
$$g_2=\frac{s}{s+s^{-1}}\raisebox{-3mm}{\epsfxsize.3in\epsffile{parra.ai}}\quad - \frac{s^{-1}}{s+s^{-1}}\quad \raisebox{-3mm}{\epsfxsize.3in\epsffile{left.ai}}.$$

Hence,
$$\raisebox{-3mm}{\epsfxsize.3in\epsffile{parra.ai}}\quad=f_2+g_2.$$

Now apply the multiplicative homomorphism $s_2: H_2\to K_2$ constructed in the proof of
\cite{BB2000}[Theorem 3.1],

$$s_2(\quad \raisebox{-3mm}{\epsfxsize.3in\epsffile{parra.ai}}\quad)=\quad\raisebox{-3mm}{\epsfxsize.3in\epsffile{parra1.ai}}-\delta^{-1} \quad\raisebox{-3mm}{\epsfxsize.3in\epsffile{parra2.ai}}\quad.$$
We obtain that,
\begin{alignat}{3}
\raisebox{-3mm}{\epsfxsize.3in\epsffile{parra1.ai}}&=s_2(\quad \raisebox{-3mm}{\epsfxsize.3in\epsffile{parra.ai}}\quad)+  \delta^{-1} \quad\raisebox{-3mm}{\epsfxsize.3in\epsffile{parra2.ai}}\notag\\
&=s_2(f_2)+s_2(g_2)+ \delta^{-1} \quad\raisebox{-3mm}{\epsfxsize.3in\epsffile{parra2.ai}}.\notag\\
&= \widetilde{y}_{\lambda}+ \widetilde{y}_{\mu}+\delta^{-1} \quad\raisebox{-3mm}{\epsfxsize.3in\epsffile{parra2.ai}}.\notag
\end{alignat}

Now we have the following skein relations,
$$\raisebox{-12mm}{\epsfxsize.0in\epsffile{ex1.ai}}=\delta\quad\raisebox{-9mm}{\epsfxsize.5in\epsffile{exlambda.ai}}+\delta\quad\raisebox{-9mm}{\epsfxsize.5in\epsffile{exmu.ai}}+\delta \delta^{-1} \quad\raisebox{-3mm}{\epsfxsize.3in\epsffile{parra2.ai}}$$

\begin{alignat}{2}
\raisebox{-12mm}{\epsfxsize.0in\epsffile{ex2.ai}}&=\quad\raisebox{-9mm}{\epsfxsize.5in\epsffile{exlambda1.ai}}+\quad\raisebox{-9mm}{\epsfxsize.5in\epsffile{exmu1.ai}}+\delta \delta^{-1} \quad\raisebox{-3mm}{\epsfxsize.3in\epsffile{parra2.ai}}\notag\\
&=(c_\lambda+\delta)\quad\raisebox{-9mm}{\epsfxsize.5in\epsffile{exlambda.ai}}+(c_\mu+\delta)\quad\raisebox{-9mm}{\epsfxsize.5in\epsffile{exmu.ai}}+\delta \delta^{-1} \quad\raisebox{-3mm}{\epsfxsize.3in\epsffile{parra2.ai}}. \notag
\end{alignat}
The last equality follows from Lemma 3.1.

From the proof of the Epimorphism Lemma 3.2, all basis elements given by the form in Theorem 5 which cut the separating disc nontrivially are $0$. Therefore modulo these basis elements of $S^1\times D^2\# S^1\times D^2$ which are zero in $S^1\times S^2\# S^1\times S^2$, then Equation (I) implies,
$$\raisebox{-7mm}{\epsfxsize.0in\epsffile{ex5.ai}}\equiv \delta^{-1} \quad\raisebox{-3mm}{\epsfxsize.3in\epsffile{parra2.ai}}$$
Thus, 
$$K_{knot}=\delta^{-2} \raisebox{-20mm}{\epsfxsize.0in\epsffile{exs1s21.ai}}=\delta^{-1}\phi.$$
i.e. $K_{knot}=\delta^{-1}\phi$ in $K(S^1\times S^2\# S^1\times S^2)$.

\section{Proof of Lemma 3.1}

\noindent {\bf Lemma 3.1:}
$$\raisebox{-15mm}{\epsfxsize.0in\epsffile{alasab.ai}}\quad=c_\lambda \quad
\raisebox{-12.5mm}{\epsfxsize.0in\epsffile{alasa.ai}}.$$

\begin{proof}
We proceed by induction on the size $|\lambda|$ of $\lambda$.

(1) When $|\lambda|=0$, it is trivial.

(2) Suppose it is true for $|\lambda|<n$. Now we assume $\lambda$ is of size $n$.

By the absorbing property given in \cite{BB2000}[Chapter 5], we have the following identity: $$\raisebox{-12mm}{\epsfxsize.0in\epsffile{alasa.ai}}\quad=\quad
\raisebox{-28mm}{\epsfxsize.0in\epsffile{alas1.ai}},$$
as well as the following skein relation \cite[Prop. 6.1]{BB2000}:
$$\raisebox{-27mm}{\ \epsfxsize0in\epsffile{alasb.ai}}
=s^{2cn(c)}\quad \raisebox{-12mm}{\
\epsfxsize0in\epsffile{alasa.ai}},$$
We develop the following relation from above,
$$\raisebox{-27mm}{\ \epsfxsize0in
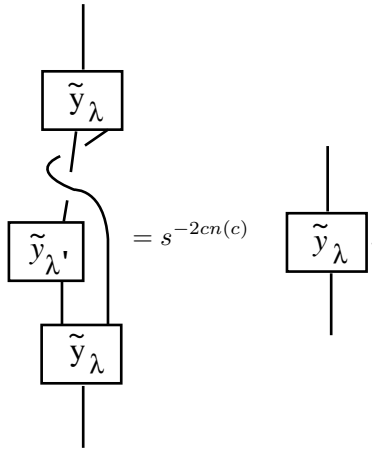
\epsffile{alasb1.ai}}
=s^{-2cn(c)}\quad \raisebox{-12mm}{\
\epsfxsize0in\epsffile{alasa.ai}}.$$
Here $cn(c)$ is the content
of the extreme cell $c$ of $\lambda$ to be removed to obtain
$\lambda'$.

Now, $$\raisebox{-28mm}{\epsfxsize.0in
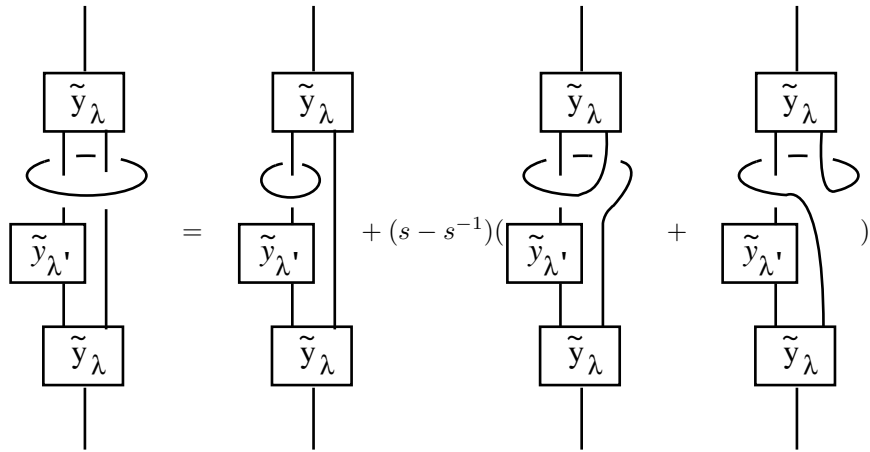
\epsffile{alas2.ai}}\quad=\quad
\raisebox{-28mm}{\epsfxsize.0in\epsffile{alas3.ai}}+(s-s^{-1})(
\raisebox{-28mm}{\epsfxsize.0in\epsffile{alas4.ai}}\quad+\quad
\raisebox{-28mm}{\epsfxsize.0in\epsffile{alas5.ai}})$$
$$=\quad
\raisebox{-28mm}{\epsfxsize.0in\epsffile{alas3.ai}}+(s-s^{-1})(\alpha
s^{2cn(c)}-\alpha^{-1}s^{-2cn(c)})
\raisebox{-12mm}{\epsfxsize.0in\epsffile{alasa.ai}}$$
(by induction on $\lambda'$)
$$=c_{\lambda'}\quad
\raisebox{-28mm}{\epsfxsize.0in\epsffile{alas1.ai}}+(s-s^{-1})(\alpha
s^{2cn(c)}-\alpha^{-1}s^{-2cn(c)})
\raisebox{-12mm}{\epsfxsize.0in\epsffile{alasa.ai}}=c_\lambda \widetilde{y}_\lambda.$$
\end{proof}

\end{document}